\documentclass{amsart}[12pt]
\usepackage{amsmath, amsthm, amscd, amsfonts,}
\usepackage{graphicx,float}
\usepackage{comment}

\usepackage{amssymb}

\usepackage{pb-diagram}

\newcommand{\url}[1]{{\tt #1}}

\DeclareMathOperator{\Col}{Col}

\def\MPB{{\mathbb{P}}}

\newtheorem{theorem}{Theorem}[section]

\newtheorem{corollary}[theorem]{Corollary}

\newtheorem{open Question}[theorem]{Open Question}

\numberwithin{equation}{section}

\def\MPB{{\mathbb{P}}}

\def\rmark{\mbox{$\rm\bf\rule{0.06em}{1.45ex}\kern-0.05em R$}}
\def\pmark{\mbox{$\rm\bf\rule{0.06em}{1.45ex}\kern-0.05em P$}}
\def\nmark{\mbox{$\rm\bf\rule{0.06em}{1.45ex}\kern-0.05em N$}}
\def\vdash{\mbox{$\rm\| \kern-0.13em -$}}

\newcommand{\GCH}{\hbox{GCH}}

\usepackage{pb-diagram}

\def\rmark{\mbox{$\rm\bf\rule{0.06em}{1.45ex}\kern-0.05em R$}}
\def\pmark{\mbox{$\rm\bf\rule{0.06em}{1.45ex}\kern-0.05em P$}}
\def\nmark{\mbox{$\rm\bf\rule{0.06em}{1.45ex}\kern-0.05em N$}}
\def\vdash{\mbox{$\rm\| \kern-0.13em -$}}

\begin{document}

\title[Notes  on countably generated complete Boolean algebras]{Notes  on countably generated complete Boolean algebras}

\author[ M. Golshani]{Mohammad Golshani }

\thanks{ The  author's research has been supported by a grant from IPM (No. 91030417).}
\maketitle

\begin{abstract}
We give a necessary and sufficient condition for an atomless Boolean algebra to be countably generated, and use it to give  new proofs of some some know facts
due to Gaifman-Hales and Solovay and also due to Jech, Kunen and Magidor. We also show that Jensen's coding theorem can be used to provide cardinal preserving countably generated
complete
Boolean algebras of arbitrary large cardinality. This answers a question of Jech \cite{jech} from 1976.
\end{abstract}

\section{Introduction}
In this paper we study atomless complete Boolean algebras (cBA for short) and provide a necessary and   sufficient condition, in terms of forcing,
for them to be countably generated. Recall that a complete Boolean algebra $\mathbb{B}$ is countably generated if there exists a countable set
$X \subseteq \mathbb{B}$ such that $X$ generates $\mathbb{B},$ i.e., $\mathbb{B}$ is the least complete subalgebra of
$\mathbb{B}$ that includes $X$. Our first result can be stated as follows.
\begin{theorem} \footnote{For application given later, just the part $(a) \to (b)$ is sufficient.}
Assume $\MPB$ is a non-trivial forcing notion and let $\mathbb{B}=R.O(\MPB)$ be the boolean completion of $\MPB.$ The following are equivalent:
\begin{enumerate}
\item [(a)] There  exists a $\mathbb{B}$-name $\dot{R}$ for a subset of $\omega$
such that for any $G$ which is $\mathbb{B}$-generic over $V$, we have $V[G]= V[\dot{R}[G]]$.

\item [(b)] $\mathbb{B}$ is countably generated.
\end{enumerate}
\end{theorem}
As an immediate corollary of the above theorem, we have the following, which gives a very simple proof of a result of Gaifman-Hales-Solovay.
\begin{corollary}
Assume $\kappa$ is an infinite cardinal and let $\mathbb{B}=R.O(\Col(\omega, \kappa)).$ Then $\mathbb{B}$ is countably generated.
\end{corollary}
As another corollary, we give a simple proof of the following result. In \cite{stavi}, this result is attributed independently to Jech, Kunen and Magidor (and possibly others).
\begin{corollary}
Assume $\kappa$ is a weakly compact cardinal. Then there is no atomless countably generated $\kappa$-c.c. complete Boolean algebra of size $\geq \kappa.$
\end{corollary}
 We may note that the countably generated complete Boolean
algebras introduced by Gaifman \cite{gaifman}, Hales \cite{hales}, Solovay \cite{solovay} and Kripke \cite{kripke} all collapse cardinals.
In \cite{jech}, Jech produced a cardinal preserving
countably generated complete Boolean algebra of size $\aleph_{\omega+2}$. He asked if we can produce such Boolean algebras of bigger cardinality.
Using (variants of) Jensen's coding
theorem, we give an affirmative answer to his question.
\begin{theorem}
Assume $\GCH$ holds. Then there are  cardinal preserving countably generated complete
Boolean algebras of arbitrary large cardinality.
\end{theorem}

\section{Characterization of countably generated cBA}
In this section we prove Theorem 1.1. and then use it to prove corollary 1.2
\begin{proof} [Proof of Theorem 1.1]
First assume  $\MPB$ is a non-trivial forcing notion and let $\mathbb{B}=R.O(\MPB)$. Also let $\dot{R}$
be a $\mathbb{B}$-name for a subset of $\omega$
such that for any $G$ which is $\mathbb{B}$-generic over $V$, we have $V[G]= V[\dot{R}[G]]$.
We show that $\mathbb{B}$ is countably generated.
Let
\[
X= \{ \parallel  \check{n} \in \dot{R} \parallel_{\mathbb{B}}~ \mid n < \omega            \}.
\]
We claim that $X$ is a set of generators. To see this, let $\mathbb{B}_X$ be the least complete subalgebra of $\mathbb{B}$
generated by $X$. We must show that $\mathbb{B}_X=\mathbb{B}$.
For any $\mathbb{B}$-generic filter $G$ over $V$ we have
\[
V[G]= = V[\dot{R}[G]] = V [G \cap \mathbb{B}_X].
\]
So by a theorem of Vop\v{e}nka \cite{vopenka}, the atomless part of $\mathbb{B}$ is isomorphic to the atomless
part of $\mathbb{B}_X.$ As $\mathbb{B}$ is atomless, it is  isomorphic to a countably
generated cBA and therefore, $\mathbb{B}$ is countably generated.

Conversely, suppose $\mathbb{B}$ is countably generated as witnessed by  a countable set $X \subseteq \mathbb{B}.$ Let $X=\{b_n \mid n < \omega   \}$
be an enumeration of $X$. Define a $\mathbb{B}$-name $\dot{R}$, for a subset of $\omega,$ such that
\begin{center}
$\Vdash_{\mathbb{B}}$``$\check{n} \in \dot{R} \Longleftrightarrow \check{b}_n \in \dot{G}$'',
\end{center}
where $\dot{G}$ is the canonical $\mathbb{B}$-name for the generic filter. Then for any $G$, which is $\mathbb{B}$-generic over $V$, we have
\[
V[G]= V[G \cap X] = V[\dot{R}[G]].
\]
So $\dot{R}$ witnesses the truth of $($a$)$.
\end{proof}

\begin{proof} [Proof of Corollary 1.2]
Let $\mathbb{B}=R.O(\Col(\omega, \kappa)).$ Also let $\dot{R}$ be a name such that
\begin{center}
$\Vdash_{\mathbb{B}}$``$\dot{R} \subseteq \omega\times \omega$ codes a well-ordering of $\omega$
of order type $\kappa$''.
\end{center}
$\dot{R}$ can be  chosen so that for all $\mathbb{B}$-generic filter $G$ over $V$, we have $V[G]=V[\dot{R}[G]]$. So by
Theorem 1.1., $\mathbb{B}=R.O(\Col(\omega, \kappa))$ is countably generated.
\end{proof}

\begin{proof} [Proof of Corollary 1.3]
Assume not, and let $\mathbb{B}$ be an atomless countably generated $\kappa$-c.c. cBA of size $\geq \kappa.$ Let $\dot{R} \subseteq \omega$
be as in Theorem 1.1. By standard facts (see \cite{hamkins}), we can find $\mathbb{Q}$, a complete subalgebra of $\mathbb{B}$, of size $<\kappa$ such that $\dot{R}$
is a $\mathbb{Q}$-name. But this is impossible, as then $\dot{R}$ can not code the generic filter for $\mathbb{B}.$
\end{proof}
\section{Cardinal preserving countably generated cBAs of arbitrary size}
In this section we prove Theorem 1.4. Our main tool will be the set version of Jensen's coding theorem.
\begin{proof} [Proof of Theorem 1.4.] Assume $\kappa$ is an infinite cardinal. We can assume that it is regular. Let $\MPB$
be a cardinal preserving  forcing notion which adds a Cohen subset of $\kappa$ and then cods it into a real in such a way that for each $\MPB$-generic filter
$G$ over $V$, we an find a real $R$ with $V[G]=V[R].$ This is possible by Jensen's coding theorem \cite{jensen} (note that we do not need to code the whole universe, so this is possible by a set forcing notion).

It is clear that $|\MPB| \geq \kappa$, and by Theorem 1.1, $\mathbb{B}=R.O.(\MPB)$ is countably generated.
\end{proof}

School of Mathematics, Institute for Research in Fundamental Sciences (IPM), P.O. Box:
19395-5746, Tehran-Iran.

E-mail address: golshani.m@gmail.com
\end{document}